\documentclass{amsart}

\usepackage{amsfonts}
\usepackage{amsmath}
\usepackage{amsthm}
\usepackage{amssymb}
\usepackage[english]{babel}
\usepackage{graphics}
\usepackage{latexsym}
\usepackage{longtable} 
\usepackage{mathrsfs} 
\usepackage{graphicx}
\usepackage{wrapfig} 
\usepackage[pdftex,colorlinks=true,linkcolor=blue,citecolor=magenta]{hyperref}
\usepackage{enumitem}
\usepackage{esint}

\newtheorem{theorem}{Theorem}
\newtheorem{corollary}[theorem]{Corollary}
\newtheorem{lemma}[theorem]{Lemma}
\newtheorem{proposition}[theorem]{Proposition}

\newtheorem{claim}[theorem]{Claim}
\newtheorem{example}[theorem]{Example}

\theoremstyle{definition}
\newtheorem{definition}[theorem]{Definition}
\newtheorem{remark}[theorem]{Remark}

\newcommand{\mH}{\mathcal{H}}

\newcommand{\mP}{\mathscr{P}}

\newcommand{\mO}{\mathcal{O}}

\newcommand{\mD}{\mathcal{D}}

\newcommand{\E}{\mathrm{E}}
\newcommand{\D}{\mathrm{D}}
\renewcommand{\H}{\mathrm{H}}

\newcommand{\A}{\mathrm{A}}
\newcommand{\R}{\mathbb{R}}

\newcommand{\mB}{\mathbb{B}}

\newcommand{\noi}{\noindent}
\newcommand{\ms}{\medskip}

\newcommand{\al}{\alpha}

\newcommand{\ga}{\gamma}
\newcommand{\Ga}{\Gamma}
\newcommand{\de}{\delta}
\newcommand{\De}{\Delta}
\newcommand{\e}{\varepsilon}
\newcommand{\si}{\sigma}

\newcommand{\la}{\lambda}

\newcommand{\Om}{\Omega}

\newcommand{\larrow}{\longrightarrow}
\newcommand{\ot}{\otimes}

\newcommand{\p}{\partial}
\newcommand{\sub}{\subseteq}

\newcommand{\by}{\times}
\newcommand{\rk}{\mathrm{rk}}

\newcommand{\diam}{\mathrm{diam}}
\newcommand{\dist}{\mathrm{dist}}

\renewcommand{\div}{\mathrm{div}}

\newcommand{\spn}{\mathrm{span}}
\newcommand{\supp}{\mathrm{supp}}

\newcommand{\bt}{\begin{theorem}}\newcommand{\et}{\end{theorem}}
\newcommand{\bd}{\begin{definition}}\newcommand{\ed}{\end{definition}}
\newcommand{\bl}{\begin{lemma}}\newcommand{\el}{\end{lemma}}
\newcommand{\beq}{\begin{equation}}\newcommand{\eeq}{\end{equation}}
\newcommand{\bc}{\begin{claim}}\newcommand{\ec}{\end{claim}}
\newcommand{\bex}{\begin{example}}\newcommand{\eex}{\end{example}}
\newcommand{\bcor}{\begin{corollary}}\newcommand{\ecor}{\end{corollary}}
\newcommand{\bp}{\begin{proof}}\newcommand{\ep}{\end{proof}}

\newcommand{\BPL}{\medskip \noindent \textbf{Proof of Lemma} }

\newcommand{\BPP}{\medskip \noindent \textbf{Proof of Proposition} }
\newcommand{\BPT}{\medskip \noindent \textbf{Proof of Theorem} }

\numberwithin{equation}{section}

\begin{document}

\title[Vectorial variational principles in $L^\infty$]{Vectorial variational principles in $L^\infty$ and their characterisation through PDE systems}

\author{Birzhan Ayanbayev}

\address{B.A.: Department of Mathematics and Statistics, University of Reading, Whiteknights, PO Box 220, Reading RG6 6AX, United Kingdom}

\email{b.ayanbayev@pgr.reading.ac.uk}

\author{Nikos Katzourakis}

\address{N.K.: (corresponding author) Department of Mathematics and Statistics, University of Reading, Whiteknights, PO Box 220, Reading RG6 6AX, United Kingdom}

\email{n.katzourakis@reading.ac.uk}

  \thanks{\!\!\!\!\!\!\!\!\texttt{N.K. has been partially financially supported by the EPSRC grant EP/N017412/1}}
  

\date{}


\begin{abstract} We discuss two distinct minimality principles for general supremal first order functionals for maps and characterise them through solvability of associated second order PDE systems. Specifically, we consider Aronsson's standard notion of absolute minimisers and the concept of $\infty$-minimal maps introduced more recently by the second author. We prove that $C^1$ absolute minimisers characterise a divergence system with parameters probability measures and that $C^2$ $\infty$-minimal maps characterise Aronsson's PDE system. Since in the scalar case these different variational concepts coincide, it follows that the non-divergence Aronsson's equation has an equivalent divergence counterpart.

\end{abstract}


\date{}


\maketitle

\section{Introduction} \label{section1}

Let $n,N\in \mathbb{N}$ and $\H \in C^2\big(\Omega\times \mathbb{R}^N\! \times\R^{N\by n}\big)$ with $\Om \subseteq \R^n$ an open set. In this paper we consider the supremal functional
\beq \label{1.1}
 \ \ \ \E_\infty (u,\mathcal{O}) := \, \underset{\mathcal{O}}{\mathrm{ess}\,\sup}\, \H (\cdot,u,\mathrm{D} u)  ,\ \ \ u\in W^{1,\infty}_\text{loc}(\Omega;\mathbb{R}^N),\ \ \mathcal{O} \Subset \Omega,
\eeq
defined on maps $u : \R^n \supseteq \Om \larrow \R^N$. In \eqref{1.1}  and subsequently, we see the gradient as a matrix map $\D u=(\D_iu_\al)_{i=1...n}^{\al=1...n} : \R^n \supseteq \Om \larrow \R^{N\by n}$. Variational problems for \eqref{1.1} have been pioneered by Aronsson in the 1960s in the scalar case $N=1$ (\cite{A1}-\cite{A5}). Nowadays the study of such functionals (and of their associated PDEs describing critical points) form a fairly well-developed area of vivid interest, called Calculus of Variations in $L^\infty$. For pedagogical general introductions to the theme we refer to \cite{ACJ,C,K4}. 

One of the main difficulties in the study of \eqref{1.1} which prevents us from utilising the standard machinery of Calculus of Variations for conventional (integral) functionals as e.g.\ in \cite{D} is that it is non-local, in the sense that a global minimisers $u$ of $\E_\infty (\cdot,\Om)$ in $W^{1,\infty}_g(\Omega;\mathbb{R}^N)$ for some fixed boundary data $g$ may not minimise $\E_\infty (\cdot,\mathcal{O})$ in $W^{1,\infty}_u(\mO;\mathbb{R}^N)$. Namely, global minimisers are not generally local minimisers, a property which is automatic for integral functionals. The remedy proposed by Aronsson (adapted) to the vector case is to build locality into the minimality notion:
\begin{definition} \label{Def1} Let $u\in W^{1,\infty}_\text{loc}(\Omega;\mathbb{R}^N)$. We say that $u$ is an {\it absolute minimiser} of \eqref{1.1} on $\Om$ if
\beq \label{1.2}
\left.
\begin{array}{l}
\forall\ \mO \Subset \Om, \\
\forall\ \phi \in W^{1,\infty}_0(\mO;\mathbb{R}^N) \\
\end{array}
\right\} \ \ \Longrightarrow \ \
\E_\infty (u,\mathcal{O})\, \leq \, \E_\infty (u+\phi,\mathcal{O}).
\eeq
\end{definition}
In the scalar case of $N=1$, Aronsson's concept of absolute minimisers turns out to be the appropriate substitute of mere minimisers. Indeed, absolute minimisers possess the desired uniqueness properties subject to boundary conditions and, most importantly, the possibility to characterise them through a necessary (and sufficient) condition of satisfaction of a certain nonlinear nondivergence second order PDE, known as the Aronsson equation (\cite{ACJS, ACJ, BEJ, B, BJW1, BJW2, CD, CDP, C1, CEG, J, MWZ, Y}). The latter can be written for functions $u\in C^2(\Om)$ as
\beq
\label{1.3}
\H_{P}(\cdot, u, \mathrm{D}u)\cdot \D\big(\H(\cdot, u, \mathrm{D} u)\big) \, = \, 0.
\eeq
The Aronsson equation, being degenerate elliptic and non-divergence when formally expanded, is typically studied in the framework of viscosity solutions. In the above, $\H_P,\H_\eta,\H_x$ denotes the derivatives of $\H(x,\eta,P)$ with respect to the respective arguments and ``$\cdot$" is the Euclidean inner product. 

In this paper we are interested in characterising appropriately defined minimisers of \eqref{1.1} in the general vectorial case of $N\geq 2$ through solvability of associated PDE systems which generalise the Aronsson equation \eqref{1.3}. As the wording suggests and we explain below, when $N\geq 2$ Aronsson's notion of Definition \ref{Def1} is no longer the unique possible $L^\infty$ variational concept. In any case, the extension of Aronsson's equation to the vectorial case reads
\beq
\label{1.4}
\begin{split}
\ \ \ &\H_{P}(\cdot, u, \mathrm{D}u)\, \D\big(\H(\cdot, u, \mathrm{D} u)\big)
\\ 
& +\, \H(\cdot, u, \mathrm{D} u) \, [\H_{P}(\cdot, u, \mathrm{D} u)]^\bot  \Big(\mathrm{Div}\big(\H_{P}(\cdot, u, \mathrm{D} u)\big)- \H_{\eta}(\cdot, u, \mathrm{D} u)\Big) =\, 0.
\end{split}
\eeq
In the above, for any linear map $A : \R^n \larrow \R^N$, $[A]^\bot$ symbolises the orthogonal projection $\mathrm{Proj}_{{\mathrm{R}}(A)^\bot}$ on the orthogonal complement of its range ${\mathrm{R}}(A) \sub \R^N$. We will refer to the PDE system \eqref{1.4} as the ``Aronsson system", in spite of the fact it was actually derived by the second author in \cite{K1}, wherein the connections between general vectorial variational problems and their associated PDEs were first studied, namely those playing the role of Euler-Lagrange equations in $L^\infty$. The Aronsson system was derived through the well-known method of $L^p$-approximations and is being studied quite systematically since its discovery, see e.g. \cite{K1}-\cite{K3}, \cite{K6, KS}. The additional normal term which is not present in the scalar case imposes an extra layer of complexity, as it might be discontinuous even for smooth solutions (see \cite{K1/2, K3}).

For simplicity and in order to illustrate the main ideas in a manner which minimises technical complications, {\it in this paper we restrict our attention exclusively to regular minimisers and solutions}. In general, solutions to \eqref{1.4} are nonsmooth and the lack of divergence structure combined with its vectorial nature renders its study beyond the reach of viscosity solutions. To this end, the theory of $\mD$-solutions introduced in \cite{K6} and subsequently utilised in several works (see e.g.\ \cite{AK,CKP,K6,K7}) offers a viable alternative for the study of general locally Lipschitz solutions to \eqref{1.4}, and in fact it works far beyond the realm of Calculus of Variations in $L^\infty$. We therefore leave the generalisation of the results herein to a lower regularity setting for future work.

Additionally to absolute minimisers, for reasons to be explained later, in the paper \cite{K2} a special case of the next $L^\infty$ variational concept was introduced (therein for $\H(x,\eta,P)=|P|^2$):

\begin{definition} 
\label{Def2} 
Let $u\in C^1(\Omega;\mathbb{R}^N)$. We say that $u$ is an {\it $\infty$-minimal map} for \eqref{1.1} on $\Om$ if (i) and (ii) below hold true:
\ms

\noi (i) $u$ is a {\it rank-one absolute minimiser}, namely it minimises with respect to essentially scalar variations vanishing on the boundary along fixed unit directions:

\beq \label{1.5}
\left.
\begin{array}{l}
\forall\ \mO \Subset \Om, \ \forall\ \xi \in \R^N\\
\forall\ \phi \in C^1_0(\overline{\mO};\spn[\xi]) \\
\end{array}
\right\} \ \Longrightarrow \ \ 
\E_\infty (u,\mathcal{O}) \leq  \E_\infty (u+\phi,\mathcal{O}).
\eeq

\ms

\noi (ii) $u$ has {\it $\infty$-minimal area}, namely it minimises with respect to variations which are normal to the range of the matrix field $\H_P(\cdot,u,\D u)$ and free on the boundary:

\beq \label{1.6}
\left.
\begin{array}{l}
\forall\ \mO \Subset \Om, \ \ \forall\ \phi \in C^1(\R^n;\R^N)\\
\text{with }\phi^\top\H_P(\cdot,u,\D u)=0 \text{ on }\mO 
\end{array}
\right\} \  \Longrightarrow \ \ 
\E_\infty (u,\mathcal{O}) \leq  \E_\infty (u+\phi,\mathcal{O}).
\eeq
\end{definition}

In the above,
\[
C^1_0(\overline{\mO};\R^N)\, :=\, \big\{\psi \in C^1(\R^n;\R^N)\ :\ \psi=0\text{ on }\p\mO \big\}.
\]
Note also that {\it when $N=1$ absolute minimisers and $\infty$-minimal maps coincide}, at least when $\{\H_P=0\}\sub \{\H=0\}$. Further, in the event that $\H_P(\cdot,u,\D u)$ has  discontinuous rank on $\mO$, the only continuous normal vector fields $\phi$ may be only those vanishing on the set of discontinuities.

In \cite{K2} it was proved that $C^2$ $\infty$-minimal maps of full rank (namely immersions or submersions) are $\infty$-Harmonic, that is solutions to the so-called $\infty$-Laplace system. The latter is a special case of \eqref{1.4}, corresponding to the choice $\H(x,\eta, P)=|P|^2$:
\beq 
\label{1.7}
\begin{split}
\ \ \ \mathrm{D}u\, \D\big(|\D u|^2\big) \, +\, |\D u|^2 \,  [ \mathrm{D} u ]^\bot \De u\, =\, 0.
\end{split}
\eeq
The fullness of rank was assumed because of the possible discontinuity of the coefficient $[\D u]^\bot$, which may well happen even for smooth solutions (for explicit examples see \cite{K1/2}). In this paper we bypass this difficulty by replacing the orthogonal projection $[\, \cdot \, ]^\bot$ by the projection on the subspace of those normal vectors which have local normal $C^1$ extensions in a open neighbourhood:
\begin{definition}
\label{def3} 
 Let $V : \R^n \supseteq \Om \larrow \R^{N\by n}$ be a matrix field and note that 
\[
{\mathrm{R}}(V(x))^\bot =\, {\mathrm{N}}(V(x)^\top), 
\]
where for any $x\in\Om$, ${\mathrm{N}}(V(x)^\top)$ is the nullspace of the transpose $V(x)^\top\in\R^{n\by N}$. We define the orthogonal projection
\[
[\![V(x)]\!]^\bot :=\, \mathrm{Proj}_{\tilde{\mathrm{N}}(V(x)^\top)}, \ \ \ [\![V(\cdot)]\!]^\bot  : \ \  \R^n \supseteq \Om \larrow \R^{N\by N},
\]
where $\tilde{\mathrm{N}}(V(x)^\top)$ is the reduced nullspace, given by
\[
\begin{split}
\tilde{\mathrm{N}}(V(x)^\top)\, :=\, \Big\{\xi \in & {\mathrm{N}}(V(x)^\top)  \ \Big| \ \exists\ \e>0 \ \, \&\ \,  \exists \ \bar \xi \in C^1(\R^n;\R^N): 
\\
& \bar\xi(x)  =\xi\ \, \&\ \, \bar \xi(y) \in {\mathrm{N}}(V(y)^\top),\ \forall\, y \in \mB_\e(x)\Big\}.
\end{split}
\]
\end{definition} 
It is a triviality to check that $\tilde{\mathrm{N}}(V(x)^\top)$ is indeed a vector space and that
\[
[\![V(x)]\!]^\bot [V(x)]^\bot=\, [\![V(x)]\!]^\bot,
\]
where $[V(x)]^\bot = \mathrm{Proj}_{{\mathrm{N}}(V(x)^\top)}$. Note that the definition could be written in a more concise manner by using the algebraic language of {\it sheaves and germs}, but we refrained from doing so as there is no real benefit in this simple case. 

The first main result in this paper is the next variational characterisation of the Aronsson system \eqref{1.4}.

\bt[Variational Structure of Aronsson's system] \label{th4}  Let $u : \R^n\supseteq \Om \larrow \R^N$ be a  map in $C^2(\Om;\R^N)$. Then:

\smallskip

\noi {\rm (I)} If $u$ is a rank-one absolute minimiser for \eqref{1.1} on $\Om$ (Definition \ref{Def2}(i)), then it solves
\beq \label{1.8}
\H_{P}(\cdot, u, \mathrm{D}u) \, \D\big(\H(\cdot, u, \mathrm{D} u)\big) \, = \, 0\ \text{ on }\Om.
\eeq
The opposite is true if in addition $\H$ does not depend on $\eta \in \R^N$ and $\H_P(\cdot,\D u)$ has full rank on $\Om$.

\ms

\noi {\rm (II)} If $u$ has $\infty$-minimal area for \eqref{1.1} on $\Om$ (Definition \ref{Def2}(ii)), then it solves
\beq \label{1.9}
\H(\cdot, u, \mathrm{D} u)\, [\![\H_{P}(\cdot, u, \mathrm{D} u)]\!]^\bot \Big(\mathrm{Div}\big(\H_{P}(\cdot, u, \mathrm{D} u)\big)- \H_{\eta}(\cdot, u, \mathrm{D} u)\Big) =\, 0\ \text{ on }\Om.
\eeq
The opposite is true if in addition for any $x\in \Om$, $\H(x,\cdot,\cdot)$ is convex on $\R^n\!\by \R^{N\by n}$.

\ms

\noi {\rm (III)} If $u$ is $\infty$-minimal map for \eqref{1.1} on $\Om$, then it solves the (reduced) Aronsson system
\[
\begin{split}
\ \ \ \A_\infty u\,:=\ &\H_{P}(\cdot, u, \mathrm{D}u)\, \D\big(\H(\cdot, u, \mathrm{D} u)\big)
\\ 
& +\, \H(\cdot, u, \mathrm{D} u) \, [\![\H_{P}(\cdot, u, \mathrm{D} u)]\!]^\bot  \Big(\mathrm{Div}\big(\H_{P}(\cdot, u, \mathrm{D} u)\big)- \H_{\eta}(\cdot, u, \mathrm{D} u)\Big) =\, 0.
\end{split}
\]
The opposite is true if in addition $\H$ does not depend on $\eta \in \R^N$, $\H_P(\cdot,\D u)$ has full rank on $\Om$ and for any $x\in\Om$ $H(x,\cdot)$ is convex in $\R^{N\by n}$.
\et

The emergence of two distinct sets of variations and a pair of separate PDE systems comprising \eqref{1.4} might seem at first glance mysterious. However, it is a manifestation of the fact that the (reduced) Aronsson system in fact consists of two linearly independent differential operators because of the perpendicularity between $[\![\H_{P}]\!]^\bot$ and $\H_P$; in fact, one may split $\A_\infty u =0$ to
\[
\left\{\ \ 
\begin{split}
\H_{P}(\cdot, u, \mathrm{D}u) \, \D\big(\H(\cdot, u, \mathrm{D} u)\big) \, &= \, 0,
\\ 
\H(\cdot, u, \mathrm{D} u)\, [\![\H_{P}(\cdot, u, \mathrm{D} u)]\!]^\bot \Big(\mathrm{Div}\big(\H_{P}(\cdot, u, \mathrm{D} u)\big)- \H_{\eta}(\cdot, u, \mathrm{D} u)\Big)\, &=\, 0.
\end{split}
\right. \ \ \
\]
Theorem \ref{th4} makes clear that Aronsson's absolute minimisers do {\bf not} characterise the Aronsson system when $N\geq 2$, at least when the additional natural assumptions hold true. This owes to the fact that, unlike the scalar case, {\it the Aronsson system admits arbitrarily smooth non-minimising solutions}, even in the model case of the $\infty$-Laplacian. For details we refer to \cite{KS}.

Since Aronsson's absolute minimisers do not characterise the Aronsson system, the natural question arises as to what is their PDE counterpart. The next theorem which is our second main result answers this question:

\bt[Divergence PDE characterisation of Absolute minimisers] \label{th5}  Let $u : \R^n\supseteq \Om \larrow \R^N$ be a  map in $C^1(\Om;\R^N)$. Fix also $\mO \Subset \Om$  and consider the following statements:

\smallskip

\noi {\rm (I)} $u$ is a vectorial minimiser of $\E_\infty(\cdot,\mO)$ in $C^1_u(\overline{\mO};\R^N)$.

\smallskip

\noi {\rm (II)} We have
\[
\max_{\mathrm{Argmax}\{\H(\cdot, u, \mathrm{D} u) \, :\, \overline{\mO}\}}\Big[\H_P(\cdot, u, \mathrm{D} u):\D \psi  \,+\, \H_{\eta}(\cdot, u, \mathrm{D} u)\cdot \psi \Big]\, \geq\, 0,
\]
for any $\psi \in C^1_0(\overline{\mO};\R^N)$.

\smallskip

\noi {\rm (III)} For any $\psi \in C^1_0(\overline{\mO};\R^N)$, there exists a non-empty compact set 
\beq \label{1.10}
\mathrm K_\psi \equiv \mathrm K\,\sub\, \mathrm{Argmax}\big\{\H(\cdot, u, \mathrm{D} u) \, :\, \overline{\mO}\big\}
\eeq
such that,
\beq  \label{1.12}
\Big(\H_P(\cdot, u, \mathrm{D} u):\D \psi  \,+\, \H_{\eta}(\cdot, u, \mathrm{D} u)\cdot \psi \Big)\Big|_{\mathrm K} =\, 0.
\eeq
Then, ${\rm (I)} \Longrightarrow {\rm (II)} \Longrightarrow {\rm (III)}$.
If additionally $\H(x,\cdot,\cdot)$ is convex on $\R^N\!\by \R^{N\by n}$ for any fixed $x\in\Om$, then ${\rm (III)} \Longrightarrow {\rm (I)}$ and all three statements are equivalent. Further, any of the statements above are deducible from the statement:

\smallskip

\noi {\rm (IV)} For any Radon probability measure $\si\in \mP(\overline{\mO})$ satisfying 
\beq \label{e5} 
\supp (\si) \,\sub \, \mathrm{Argmax}\big\{\H(\cdot, u, \mathrm{D} u) \, :\, \overline{\mO}\big\},
\eeq
we have 
\beq \label{1.11}
-\div\big(\H_P(\cdot, u, \mathrm{D} u)\si\big) \,+\, \H_{\eta}(\cdot, u, \mathrm{D} u)\si =\, 0,
\eeq
in the dual space $(C^1_0(\overline{\mO};\R^N))^*$.

\smallskip

\noi Finally, all statement are equivalent if $\mathrm K = \mathrm{Argmax}\big\{\H(\cdot, u, \mathrm{D} u) \, :\, \overline{\mO}\big\}$ in {\rm (III)} (this happens for instance when the argmax is a singleton set). 
\et

The result above provides an interesting characterisation of Aronsson's concept of Absolute minimisers in terms of divergence PDE systems with measures as parameters. The exact distributional meaning of \eqref{1.11} is
\[
 \int_{\overline{\mO}}\Big(\H_P(\cdot, u, \mathrm{D} u):\D \psi  \,+\, \H_{\eta}(\cdot, u, \mathrm{D} u)\cdot \psi \Big) \mathrm{d}\si =\, 0
\]
for all $\psi \in C^1_0(\overline{\mO};\R^N)$, where the ``$:$" notation in the PDE symbolises the Euclidean (Frobenius) inner product in $\R^{N\by n}$.

The idea of Theorem \ref{th5} is inspired by the paper \cite{EY} of Evans and Yu, wherein a particular case of the divergence system is derived (in the special scalar case $N=1$ for the $\infty$-Laplacian and only for $\Om=\mO$), as well as by new developments on higher order Calculus of variations in $L^\infty$ in \cite{KM,KP,PP}. 

Note that, it does not suffice to consider only $\Om=\mO$ as in \cite{EY} in order to describe absolute minimisers. For a subdomain $\mO \sub \Om$, it may well happen that the only measure $\si$ ``charging" the points of $\overline{\mO}$ where the energy density $\H(\cdot, u, \mathrm{D} u)$ is maximised is the Dirac measure at a single point $x \in \p\mO$. This is for instance the case for the standard ``Aronsson solution" of the $\infty$-Laplacian on $\R^2$, given by $u(x,y)=|x|^{4/3}-|y|^{4/3}$, as well as for any other $\infty$-Harmonic function which is nowhere Eikonal (i.e.\ $|\D u|$ is non-constant on all open subsets).


%
%







We conclude this introduction by noting that the two vectorial variational concepts we are considering herein (Definitions \ref{Def1}-\ref{Def2}) do not exhaust the plethora variational concepts in $L^\infty$. In particular, in the paper \cite{SS} the concept of {\it tight maps} was introduced in the case of $\H(x,\eta,P)=\|P\|$ where $\|\, \cdot\,\|$ is the operator norm on $\R^{N\by n}$. Additionally, in the papers \cite{AK,K7} a concept of special affine variations was considered which also characterises the Aronsson system, in fact in the generality of merely locally Lipschitz $\mD$-solutions. Finally, in the paper \cite{AB} new concepts of absolute minimisers for constrained minimisation problems have been proposed, whilst results relevant to variational principles in $L^\infty$ and applications appear in \cite{BN, BP, CDP, GNP, P, RZ}.

\section{Proofs and a maximum-minimum principle for $\H(\cdot,u,\D u)$} \label{section2}

In this section we prove our main results Theorems \ref{th4}-\ref{th5}. Before delving into that, we establish a result of independent interest, which generalises a corresponding result from \cite{K2}.

\begin{proposition}[Maximum-Minimum Principles] Suppose \label{pr1} Let $u \in C^2(\Om;\R^N)$ be a solution to \eqref{1.8}, such that $\H$ satisfies 

\noi  {\rm (a)} $\H_P(\cdot,u,\D u)$ has full rank on $\Om$,

\noi  {\rm (b)} there exists $c>0$ such that 
\[
\big(\xi^\top \H_P(x,\eta,P)\big) \cdot \big(\xi^\top P)\, \geq c \big|\xi^\top \H_P(x,\eta,P)\big|^2,
\]
for all $\xi \in \R^N$ and all $(x,\eta,P)\in \Om \by \R^N \!\by \R^{N\by n}$.

\smallskip 
 
Then, for any $\mO \Subset \Om$ we have:
 \begin{align} 
\sup_{\mO}\H(\cdot,u,\D u)\, &= \, \max_{\p \mO}\H(\cdot,u,\D u), \label{2.1}\\
\inf_{\mO}\H(\cdot,u,\D u)\, &= \, \min_{\p \mO}\H(\cdot,u,\D u). \label{2.2}
\end{align}
\end{proposition}

The proof is based on the usage of the following flow with parameters:
\begin{lemma} \label{l5}
Let $u \in C^2(\Om;\R^N)$. Consider the parametric ODE system
\beq \label{2.3}
\left\{
\begin{array}{l}
\dot{\ga}(t)\ = \ \xi^\top \H_P(\cdot,u,\D u)\big|_{\ga(t)}, \ \ t\neq 0,\ms\\
\ga(0)\ = \ x,
\end{array}
\right.
\eeq
for given $x\in \Om$ and $\xi \in \R^N$. Then, we have
\begin{align}
\ \ \ \ \ \frac{d}{d t} \Big(\H(\cdot,u,\D u)\big|_{\ga(t)}\Big)\, &= \, \xi^\top 
\H_{P}(\cdot, u, \mathrm{D}u) \, \D\big(\H(\cdot, u, \mathrm{D} u)\big)\big|_{\ga(t)} , \label{2.4}\\
\frac{d}{d t}\xi^\top u\big(\ga(t)\big)\, & \geq \, c\Big|\xi^\top \H_P(\cdot,u,\D u)\big|_{\ga(t)}\Big|^2. \label{2.5}
\end{align}
\end{lemma}

\BPL \ref{l5}. The identity \eqref{2.4} follows by a direct computation and \eqref{2.3}. For the inequality \eqref{2.5}, we have 
\[
\begin{split}
\frac{d}{d t}\xi^\top u\big(\ga(t)\big)\, &=\,  \Big(\xi^\top \D u\big(\ga(t)\big)\Big) \cdot \dot \ga(t)
\\
& = \, \Big(\xi^\top \D u\big(\ga(t)\big)\Big) \cdot \Big(\xi^\top \H_P(\cdot,u,\D u)\big|_{\ga(t)}\Big)
\\
& \geq \, c\Big|\xi^\top \H_P(\cdot,u,\D u)\big|_{\ga(t)}\Big|^2.
\end{split}
\]
The lemma ensues. \qed
\ms

\BPP \ref{pr1}. Fix $\mO \Subset \Om$. Without loss of generality, we may suppose $\mO$ is connected. Consider first the case where $\rk\big(\H_P(\cdot,u,\D u)\big) \equiv n\leq N$. Then, the matrix-valued map $\H_P(\cdot,u,\D u)$ is pointwise left invertible. Therefore, by \eqref{1.8},
\[
\big(\H_{P}(\cdot, u, \mathrm{D}u)\big)^{-1}\H_{P}(\cdot, u, \mathrm{D}u) \, \D\big(\H(\cdot, u, \mathrm{D} u)\big)\, =\, 0
\]
which, by the connectivity of $\mO$, gives $\H(\cdot, u, \mathrm{D} u)  \equiv  \mathrm{const}$ on $\mO$. The latter equality readily implies the desired conclusion. Consider now the case where $\rk\big(\H_P(\cdot,u,\D u)\big) \equiv N\leq n$. Fix $x\in \mO$ and a unit vector $\xi \in \R^n$ and consider the parametric ODE system \eqref{2.3} of Lemma \ref{l5}. By the fullness of the rank of $\H_P(\cdot,u,\D u)\big)$, we have that
\[
\big| \xi^\top \H_P(\cdot,u,\D u)\big)\big| \, \geq\, c_1 >0\ \ \text{ on }\mO.
\]
We will now show that the trajectory $\ga(t)$ reaches $\p \mO$ in finite time. To this end, we estimate
\[
\begin{split}
\|\D u\|_{L^\infty(\mO)}\diam(\mO)\, & \geq \, \|\D u\|_{L^\infty(\mO)} \big|\ga(t)\, - \, \ga(0)\Big|
\, \geq\, \bigg|\frac{\mathrm d}{\mathrm d t}\Big|_{\hat{t}}\, \xi^\top u( \ga(t))\bigg|\, t,
\end{split}
\]
for some $\hat t \in (0,t)$, by the mean value theorem. Hence,
\[
\begin{split}
\|\D u\|_{L^\infty(\mO)}\diam(\mO)\, & \geq\, \bigg|\frac{\mathrm d}{\mathrm d t}\Big|_{\hat{t}}\, \xi^\top u( \ga(t))\bigg|\, t
\\
& =\, \Big|\xi^\top \D u( \ga(\hat{t})) \cdot \dot \ga(\hat{t})\Big|\, t
\\
& =\, \Big|\xi^\top \D u( \ga(\hat{t})) \cdot \Big(\xi^\top\H_P(\cdot,u,\D u)\big|_{\ga(\hat{t})} \Big) \Big|\, t
\\
&\geq \, c_0 \Big|\xi^\top\H_P(\cdot,u,\D u)\big|_{\ga(\hat{t})}\Big|^2 t
\\
&\, \geq (c_0c_1^2)\, t.
\end{split}
\]
This proves the desired claim. Further, since $u$ solves \eqref{1.8}, by \eqref{2.4} of Lemma \ref{l5} it follows that $\H(\cdot,u,\D u)$ is constant along the trajectory. Thus, if $x \in \mO$ is chosen as a point realising either the maximum or the minimum in $\overline \mO$, then by moving along the trajectory, we reach a point $y\in\p\mO$ such that $\H(\cdot,u,\D u)\big|_x = \H(\cdot,u,\D u)\big|_y$. This establishes both the maximum and minimum principle. The proposition ensues.                      
\qed

\begin{remark}[Danskin's theorem] \label{Danskin's theorem} The central ingredient in the proofs of Theorems \ref{th4}-\ref{th5} is the next consequence of Danskin's theorem: for any $\mO \Subset \Om$ and any $u,\phi \in C^1(\Om;\R^N)$, we have the identities
\beq
\label{e0}
\left\{
\begin{split}
\frac{\mathrm d}{\mathrm d t}\Big|_{t=0^+} \! \E_\infty(u+t\phi,\mO)\, &=\, \max_{\mO(u)}\Big( \H_P(\cdot,u,\D u) : \D \phi\,+\,  \H_\eta(\cdot,u,\D u) \cdot \phi \Big),
\\
\frac{\mathrm d}{\mathrm d t}\Big|_{t=0^-} \! \E_\infty(u+t\phi,\mO)\, &=\, \min_{\mO(u)}\Big( \H_P(\cdot,u,\D u) : \D \phi\,+\,  \H_\eta(\cdot,u,\D u) \cdot \phi \Big),
\end{split}
\right.
\eeq
where
\[
\mO(u)\,:=\, \mathrm{Argmax}\big\{\H(\cdot,u,\D u)\,:\, \overline \mO \big\}.
\]
Indeed, by \cite[Theorem 1, page 643]{D} and the chain rule we have
\[
\begin{split}
\frac{\mathrm d}{\mathrm d t}\Big|_{t=0^+} \E_\infty(u+t\phi,\mO)\, &=\, \frac{\mathrm d}{\mathrm d t}\Big|_{t=0^+} \Big(\max_{\overline \mO} \H\big(\cdot,u+t\phi,\D u+t\D \phi\big)\Big)
\\
&=\, \max_{\mO(u)}\bigg(\frac{\mathrm d}{\mathrm d t}\Big|_{t=0^+} \H\big(\cdot,u+t\phi,\D u+t\D \phi\big)\bigg)
\\
& =\, \max_{\mO(u)}\Big( \H_P(\cdot,u,\D u) : \D \phi\,+\,  \H_\eta(\cdot,u,\D u) \cdot \phi \Big).
\end{split}
\]
This establishes the first identity of \eqref{e0}. The second one follows through the substitutions $\phi \leadsto -\phi$, $t \leadsto -t$. 
\end{remark}

Now we may establish Theorem \ref{th4}.

\BPT \ref{th4}. (I) Suppose first that $u$ is a rank-one absolute minimiser on $\Om$. The aim is to show that \eqref{1.8} is satisfied on $\Om$. This conclusion in fact follows by the results in \cite{K1}, but below we provide a new shorter proof. To this end, fix $x\in \Om$ and $\rho \in (0,\dist(x,\p\Om))$ and let $\mO:=\mB_\rho (x)$. We fix also $\xi\in \R^N$ and choose
\[
\phi(y)\,:=\, \xi\big( |y-x|^2-\rho^2\big).
\]
Then, $\phi \in C^1_0\big(\bar{\mB}_\rho(x);\spn[\xi]\big)$. By Remark \ref{Danskin's theorem} and our minimality assumption, the definition of one-sided derivatives yields
\beq \label{e2a}
\frac{\mathrm d}{\mathrm d t}\Big|_{t=0^-} \!\E_\infty(u+t\phi,\mO) \,\leq\, 0 \, \leq\, \frac{\mathrm d}{\mathrm d t}\Big|_{t=0^+} \! \E_\infty(u+t\phi,\mO).
\eeq
Hence, by \eqref{e2a}, \eqref{e0} and continuity there exists a point ${x_\rho}$ with $|x_\rho -x|\leq \rho$ which lies in the argmax set 
\[
(\mB_\rho(x))(u) \, =\, \mathrm{Argmax}\big\{\H(\cdot,u,\D u)\, :\, \bar \mB_\rho(x)\big\}
\]
such that 
\beq \label{e4}
\Big(\H_P(\cdot,u,\D u) : \D \phi \, +\, \H_\eta(\cdot,u,\D u) \cdot \phi\Big)\Big|_{x_\rho}\, =\,0.
\eeq
Therefore,
\beq \label{2.8a}
\xi^\top \Big(2\H_P(\cdot,u,\D u)\big|_{x_\rho} (x_\rho -x) \,+\,  \H_\eta(\cdot,u,\D u)\big|_{x_\rho}\big( |x_\rho-x|^2-\rho^2\big) \Big) \, =\, 0.
\eeq
If $x_\rho$ lies in the interior of $\mB_\rho(x)$, then it is an interior maximum and therefore 
\[
\D\big(\H(\cdot,u,\D u)\big)\big|_{x_\rho} =\, 0.
\]
This means that \eqref{1.8} is satisfied at $x_\rho$. If $x_\rho$ lies on the boundary of $\mB_\rho(x)$, then this means that
\[
\forall \ y \in \bar \mB_\rho(x), \text{ we have } \H(\cdot,u,\D u)\big|_{y} \leq\, \H(\cdot,u,\D u)\big|_{x_\rho}.
\]
The above can be rewritten as
\[
\bar \mB_\rho(x) \, \sub \, \mH(x_\rho)\,:=\, \Big\{ \H(\cdot,u,\D u) \leq\, \H(\cdot,u,\D u)\big|_{x_\rho}\Big\},
\]
and note also that $x_\rho \in \p \mB_\rho(x) \cap \p \mH(x_\rho)$. Hence, the sublevel set $\mH(x_\rho)$ satisfied an interior sphere condition at $x_\rho$. If $\D\big(\H(\cdot,u,\D u)\big)\big|_{x_\rho} =0$ then \eqref{1.8} is again satisfied at $x_\rho$. If on the other hand
\[
\D\big(\H(\cdot,u,\D u)\big)\big|_{x_\rho} \neq \, 0
\]
then $\p \mH(x_\rho)$ is a $C^1$ manifold near $x_\rho$ and the gradient above is the normal vector at the point $x_\rho$. Due to the interior sphere condition, this implies that this is also the normal vector to the sphere $\p \mB_\rho(x)$ at $x_\rho$. Thus, there exists $\la\neq 0$ such that
\beq \label{2.9a}
x_\rho -x\, =\, \la\D\big(\H(\cdot,u,\D u)\big)\big|_{x_\rho}.
\eeq
By inserting \eqref{2.9a} into \eqref{2.8a} and noting that $|x_\rho -x|=\rho$, we infer that
\[
2\la\, \xi^\top \Big(\H_P(\cdot,u,\D u) \D\big(\H(\cdot,u,\D u)\big) \Big)\big|_{x_\rho}  =\, 0.
\]
By dividing by $2\la$ and letting $\rho \to 0$, we deduce that \eqref{1.8} is satisfied at the arbitrary $x\in \Om$.

\smallskip

Conversely, suppose that $u$ satisfies \eqref{1.8} on $\Om$, together with the additional assumptions of the statement. Fix $\mO \Subset \Om$ and $\phi \in C^1_0(\overline{\mO};\spn[\xi])$. Without loss of generality, we may suppose $\mO$ is connected. Since $\phi = (\xi^\top \phi)\xi$, for convenience we set $g:=\xi^\top \phi$ and then we may write $\phi=g \xi $ with $g\in C^1_0(\overline{\mO})$. Then, the matrix-valued map $\H_P(\cdot,\D u)$ is pointwise left invertible. Therefore, by \eqref{1.8}
\[
\big(\H_{P}(\cdot, \mathrm{D}u)\big)^{-1}\H_{P}(\cdot, \mathrm{D}u) \, \D\big(\H(\cdot, \mathrm{D} u)\big) =\, 0\ \text{ on }\mO,
\]
which, by the connectivity of $\mO$, gives
\[
\H(\cdot, \mathrm{D} u)  \equiv \, \mathrm{const}\ \text{ on } \mO.
\]
Since $g\in C^1(\R^n)$ with $g=0$ on $\p\mO$, there exists at least one interior critical point $\bar x \in \mO$ such that $\D g(\bar x)=0$. By the previous, we have
\[
\begin{split}
\E_\infty(u,\mO)\, &=\, \H\big(\bar x, \D u(\bar x)\big)
\\
&=\, \H\Big(\bar x, \D u(\bar x)+\xi \ot \D g(\bar x)\Big)
\\
&=\, \H\Big(\bar x, \D u(\bar x)+\D \phi(\bar x)\Big)
\\
&\leq \, \sup_{x\in \mO}\H\Big(x, \D u(x)+\D \phi(x)\Big)
\\
&=\, \E_\infty(u+\phi,\mO).
\end{split}
\]
The conclusion ensues.
\ms

\noi (II) Suppose that $u$ has $\infty$-minimal area. Fix $x\in \Om$ and $\rho \in (0,\dist(x,\p\Om))$. Fix 
\[
\xi \, \in \, \tilde{\mathrm{N}}\Big(\H_P(\cdot,u,\D u)^\top \big|_x\Big),
\]
noting also that by Definition \ref{def3} the above set is the reduced nullspace of $\H_P(\cdot,u,\D u)^\top$ at $x$. This implies that there exists a $C^1$ extension $\bar \xi \in C^1(\R^n;\R^N)$ such that $\bar \xi(x)=\xi$ and $(\bar \xi)^\top \H_P(\cdot,u,\D u)=0$ on the closed ball $\bar \mB_\e(x)$ for some $\e \in (0,\rho)$. By differentiating the relation $(\bar \xi)^\top \H_P(\cdot,u,\D u)=0$ and taking its trace, we obtain
\beq \label{2.8}
\bar \xi \cdot \div\big(\H_P(\cdot,u,\D u)\big) \,+\, \D\bar \xi: \H_P(\cdot,u,\D u)\, =\, 0,
\eeq
on $\bar \mB_\e(x)$. Since $u$ has $\infty$-minimal area and $\bar \xi$ is an admissible normal variation, by using Remark \ref{Danskin's theorem} and arguing as in the beginning of part (I), it follows that
\beq  \label{2.9}
\Big(\bar \xi \cdot \H_\eta (\cdot,u,\D u) \,+\, \D\bar \xi: \H_P(\cdot,u,\D u)\Big)\Big|_{x_\e} =\, 0
\eeq
for some $x_\e \in (\mB_\e(x))(u)$, where 
\[
(\mB_\e(x))(u)\, =\, \mathrm{Argmax}\big\{\H(\cdot,u,\D u) \,:\, \bar \mB_\e(x)\big\}.
\]
By \eqref{2.8}-\eqref{2.9}, we infer that
\[
\bar \xi(x_\e) \cdot \Big(\div\big(\H_P(\cdot,u,\D u)\big) \,-\, \H_\eta (\cdot,u,\D u)\Big)\Big|_{x_\e} =\, 0
\]
and by letting $\e\to 0$, we deduce that
\[
\xi \cdot \Big(\div\big(\H_P(\cdot,u,\D u)\big) \,-\, \H_\eta (\cdot,u,\D u)\Big)\Big|_x\, =\, 0,
\]
for any $\xi \in \tilde{\mathrm{N}}\big(\H_P(\cdot,u,\D u)^\top \big|_x\big)$. Hence, $u$ satisfies \eqref{1.9} at the arbitrary $x\in \Om$.

\smallskip

Conversely, suppose that $u$ solves \eqref{1.9} on $\Om$. Fix $\mO\Subset \Om$ and $\phi \in C^1(\R^n;\R^N)$ such that $\phi^\top \H_P(\cdot,u,\D u)=0$ on $\mO$. Note further that by the continuity up to the boundary of all functions involved, the latter identity in fact holds on $\overline \mO$. By the satisfaction of \eqref{1.9} and Definition \ref{def3}, it follows that
\[
\phi \cdot \Big(\div\big(\H_P(\cdot,u,\D u)\big) \,-\, \H_\eta (\cdot,u,\D u)\Big)\, =\, 0,
\]
on $\overline \mO \sub \Om$. By differentiating $\phi^\top \H_P(\cdot,u,\D u)=0$, we obtain
\[
\phi \cdot \div\big(\H_P(\cdot,u,\D u)\big) \,+\, \D\phi: \H_P(\cdot,u,\D u)\, =\, 0,
\]
on $\overline \mO$. By the above two identities, we deduce
\[
\phi \cdot \H_\eta (\cdot,u,\D u) \,+\, \D\phi: \H_P(\cdot,u,\D u)\, = \, 0,
\]
on $\overline \mO$. Since $\mO(u) \sub \overline \mO$, Remark \ref{Danskin's theorem} yields that $u$ is a critical point since the left and right derivative of $\E_\infty(u+t\phi,\mO)$ at $t=0$ coincide and vanish. Since by assumption $\H(x,\cdot,\cdot)$ is convex on $\R^N\!\by \R^{N\by n}$, it follows that $\E_\infty(\cdot,\mO)$ is convex on $C^1(\overline{\mO};\R^N)$. Hence, the critical point $u$ is in fact a minimum point for this class of variations. This establishes our claim.

\ms

\noi (III) This is an immediate corollary of items (I) and (II).
\qed
\ms

Now we conclude by establishing Theorem \ref{th5}.

\BPT \ref{th5}. Fix $\mO \Subset \Om$ and $u,\phi \in C^1(\Om;\R^N)$. We show that (I) $\Longrightarrow$ (II) $\Longrightarrow$ (III) and that (III) $\Longrightarrow$ (I) under the additional convexity assumption. By recalling Remark \ref{Danskin's theorem}, note that if 
\beq \label{e1}
\E_\infty(u+t\phi,\mO) \, \geq\, \E_\infty(u,\mO) , \ \text{ for all $t\in\R$},
\eeq
then directly by \eqref{e1} and the definition of one-sided derivatives, we have
\beq \label{e2}
\frac{\mathrm d}{\mathrm d t}\Big|_{t=0^-} \!\E_\infty(u+t\phi,\mO) \,\leq\, 0 \, \leq\, \frac{\mathrm d}{\mathrm d t}\Big|_{t=0^+} \! \E_\infty(u+t\phi,\mO).
\eeq
This shows (I) $\Longrightarrow$ (II). If (II) holds, note that one also has that 
\[
\min_{\mathrm{Argmax}\{\H(\cdot, u, \mathrm{D} u) \, :\, \overline{\mO}\}}\Big[\H_P(\cdot, u, \mathrm{D} u):\D \phi  \,+\, \H_{\eta}(\cdot, u, \mathrm{D} u)\cdot \phi \Big]\, \leq\, 0,
\]
for any $\phi \in C^1_0(\overline{\mO};\R^N)$. By \eqref{e0} we see that \eqref{e2} is satisfied and by continuity we obtain the existence of a non-empty compact set $\mathrm K=\mathrm K_\phi \sub \mO(u)$ such that 
\beq \label{e4}
\Big(\H_P(\cdot,u,\D u) : \D \phi \, +\, \H_\eta(\cdot,u,\D u) \cdot \phi\Big)\Big|_{\mathrm K}\, =\,0.
\eeq
Hence, (III) ensues. If now \eqref{e4} holds true for some non-empty compact set $\mathrm K \sub \mO(u)$, then by \eqref{e0} we have that \eqref{e2} is true. If further $\H(x,\cdot,\cdot)$ is convex for all $x\in \Om$, then by Lemma \ref{Lemma-jets} given right after the proof, $t\mapsto \E_\infty(u+t\phi,\mO)$ is minimised at $t=0$ and \eqref{e1} holds true.

\smallskip

\noi (IV) $\Longrightarrow$ (III): Let $\si \in \mP(\overline \mO)$ be any Radon probability measure satisfying \eqref{e5}. Then, by assumption
\[
 \int_{\overline{\mO}}\Big(\H_P(\cdot, u, \mathrm{D} u):\D \phi  \,+\, \H_{\eta}(\cdot, u, \mathrm{D} u)\cdot \phi \Big) \mathrm{d}\si =\, 0
\]
for all $\phi \in C^1_0(\overline{\mO};\R^N)$. Fix any point $\bar x \in \mO(u)$. By choosing the Dirac measure $\bar \si \in \mP(\overline \mO)$ given by
\[
\bar \si\,:=\, \de_{\bar x}
\]
which evidently satisfies $\supp(\bar \si) =\{\bar x\} \sub \mO(u)$, we obtain
\[
\begin{split}
 \Big(\H_P(\cdot,u,\D u) : \D \phi\, &+\,  \H_\eta(\cdot,u,\D u) \cdot \phi \Big)\Big|_{\bar x}
\\
& =  \int_{\overline{\mO}}\Big(\H_P(\cdot, u, \mathrm{D} u):\D \phi  \,+\, \H_{\eta}(\cdot, u, \mathrm{D} u)\cdot \phi \Big) \mathrm{d}\bar \si 
\\
& =\, 0,
\end{split}
\]
for any $\bar x \in \mO(u)$. The conclusion ensues with $\mathrm K = \mO(u)$.

\smallskip

\noi (III) $\Longrightarrow$ (IV): If we have $\mathrm K = \mO(u)$ and
\[
\Big(\H_P(\cdot, u, \mathrm{D} u):\D \phi  \,+\, \H_{\eta}(\cdot, u, \mathrm{D} u)\cdot \phi\Big)\Big|_{\mathrm K} =\, 0,
\]
then for any Radon probability measure $\si \in \mP(\overline \mO)$ with $\supp (\si) \sub \mathrm K$, we have
\[
 \int_{\overline{\mO}}\Big(\H_P(\cdot, u, \mathrm{D} u):\D \phi  \,+\, \H_{\eta}(\cdot, u, \mathrm{D} u)\cdot \phi \Big) \mathrm{d}\si =\, 0
\]
for all $\phi \in C^1_0(\overline{\mO};\R^N)$. Hence, we have shown that
\[
-\div\big(\H_P(\cdot, u, \mathrm{D} u)\si\big) \,+\, \H_{\eta}(\cdot, u, \mathrm{D} u)\si =\, 0,
\]
in the dual space $(C^1_0(\overline{\mO};\R^N))^*$.
\qed
\ms

The next result which was utilised in the proof of Theorem \ref{th5} completes our arguments.

\begin{lemma} \label{Lemma-jets} Let $f: \R \larrow \R$ be a convex function. If the one-sided derivatives $f'(0^\pm)$ exist and $f'(0^-)\leq 0 \leq f'(0^+)$, then $f(0)$ is the global minimum of $f$ on $\R$.
\end{lemma}

\BPL \ref{Lemma-jets}. By the convexity of $f$ on $\R$, for any fixed $s\in\R$ there exists a sub-differential $p_s \in\R$ such that
\beq \label{e3}
f(t) - f(s)\, \geq \, p_s(t-s), \ \text{ for all }t \in\R.
\eeq
For the choice $t=0$ and $s>0$, we have
\[
\frac{f(s) - f(0)}{s}\, \leq \, p_s
\]
and note also that since convex functions are locally Lipschitz, the set  $(p_{s})_{0<s<1}$ is bounded. Thus, since $f'(0^+)$ exists and is non-negative, the above inequality yields
\[
0 \, \leq \, f'(0^+) \, \leq\, \liminf_{s\to 0^+}p_s \, < \, \infty.
\]
Hence, by passing to the limit as $s \to 0^+$ in the inequality \eqref{e3} for $t>0$ fixed, we obtain $f(t)-f(0)\geq 0$. The case of $t<0$ follows by arguing similarly.
\qed
\ms

\noi {\bf Conflict of Interest}: The authors declare that they have no conflict of interest.

\ms

\noi \textbf{Acknowledgement.} N.K.\ would like to thank Roger Moser, Giles Shaw and Tristan Pryer for their inspiring scientific discussion about $L^\infty$ variational problems.

\ms

\bibliographystyle{amsplain}

\begin{thebibliography}{30}

\bibitem{AK} B. Ayanbayev, N. Katzourakis, \emph{A Pointwise Characterisation of the PDE system of vectorial Calculus of variations in $L^\infty$}, Proc.  Royal Soc. Edinburgh A, in press.


\bibitem{A1} G. Aronsson, \emph{Minimization problems for the functional $sup_x F(x,
f(x), f'(x))$}, Arkiv f\"ur Mat. 6 (1965), 33 - 53.

\bibitem{A2} G. Aronsson, \emph{Minimization problems for the functional $sup_x F(x,
f(x), f'(x))$ II}, Arkiv f\"ur Mat. 6 (1966), 409 - 431.

\bibitem{A3} G. Aronsson, \emph{Extension of functions satisfying Lipschitz conditions}, Arkiv f\"ur Mat. 6 (1967), 551 - 561.

\bibitem{A4} G. Aronsson, \emph{On the partial differential equation $u_x^2 u_{xx} + 2u_x u_y u_{xy} + u_y^2 u_{yy} = 0$}, Arkiv f\"ur Mat. 7
(1968), 395 - 425.

\bibitem{A5} G. Aronsson, \emph{Minimization problems for the functional $sup_x F(x,
f(x), f'(x))$ III}, Arkiv f\"ur Mat. (1969), 509 - 512.

\bibitem{A} G. Aronsson, \emph{On Certain Minimax Problems and Pontryagin's Maximum Principle}, Calculus of Variations and PDE 37, 99 - 109 (2010).

\bibitem{AB} G. Aronsson, E.N. Barron, \emph{$L^\infty$ Variational Problems with Running Costs and Constraints}, Appl. Math. Optim. 65, 53 - 90 (2012).

\bibitem{ACJS} S.N. Armstrong, M.G. Crandall, V. Julin, C.K. Smart, \emph{Convexity Criteria and Uniqueness of Absolutely Minimising Functions}, Arch. Rational Mech. Anal. 200, 405--443 (2011). 

\bibitem{ACJ} G. Aronsson, M. Crandall, P. Juutinen \emph{A tour of the theory of
absolutely minimizing functions}, Bulletin of the AMS, New Series 41, 439--505 (2004).

\bibitem{B} E.N. Barron, \emph{Viscosity solutions and analysis in $L^\infty$}. In: Clarke F.H., Stern R.J., Sabidussi G. (eds) Nonlinear Analysis, Differential Equations and Control. NATO Science Series (Series C: Mathematical and Physical Sciences), vol 528. Springer, Dordrecht, 1999.

\bibitem{BEJ} E. N. Barron, L. C. Evans, R. Jensen, \emph{The infinity Laplacian, Aronsson's equation and their generalizations}, Trans. Amer. Math. Soc. 360 (2008), 77 - 101.

\bibitem{BJW1} E. N. Barron, R. Jensen and C. Wang, \emph{The Euler equation and absolute minimisers of $L^\infty$ functionals}, Arch. Rational Mech. Analysis 157, 255--283 (2001).

\bibitem{BJW2} N. Barron, R. Jensen, C. Wang, \emph{Lower Semicontinuity of $L^\infty$ Functionals}, Ann. I. H. Poincar\'e 18, 495--517  (2001).

\bibitem{BN} M. Bocea, V. Nesi, \emph{$\Gamma$-convergence of power-law functionals, variational principles in $L^\infty$ and applications}, SIAM J. Math. Anal. 39, 1550Ð1576 (2008).

\bibitem{BP} M. Bocea, C. Popovici, \emph{Variational principles in $L^\infty$ with applications to antiplane shear and plane stress plasticity}, Journal of Convex Analysis Vol. 18 No. 2, (2011) 403-416.

\bibitem{CD} T. Champion, L. De Pascale, \emph{Principles of comparison with distance functions for absolute minimizers}, J. Convex Anal. 14, 515Ð541 (2007).

\bibitem{CDP} T. Champion, L. De Pascale, F. Prinari, \emph{$\Ga$-convergence and absolute minimizers for supremal functionals}, COCV ESAIM: Control, Optimisation and Calculus of Variations (2004), Vol. 10, 14Ð27.

\bibitem{C} M. G. Crandall, \emph{A visit with the $\infty$-Laplacian}, in \emph{Calculus of Variations and Non-Linear Partial Differential Equations}, Springer Lecture notes in Mathematics 1927, CIME, Cetraro Italy 2005.

\bibitem{C1} M. Crandall, \emph{An Efficient Derivation of the Aronsson Equation}, Archive for Rational Mechanics and Analysis 167 (4), 271Ð279 (2003).

\bibitem{CEG} M.G. Crandall, L.C. Evans, R. Gariepy, \emph{Optimal Lipschitz extensions and the infinity Laplacian}, Calc. Var. PDE 13, 123Ð139 (2001).

\bibitem{CKP} G. Croce N. Katzourakis, G. Pisante, \emph{$\mathcal{D}$-solutions to the system of vectorial Calculus of Variations in $L^\infty$ via the singular value problem}, Discrete and Continuous Dynamical Systems 37:12, 6165-6181 (2017).

\bibitem{D} B. Dacorogna,  \emph{Direct Methods in the Calculus of Variations}, $2$nd Edition, Volume 78, Applied Mathematical Sciences, Springer, 2008.

\bibitem{D} J.M. Danskin, \emph{The theory of min-max with application}, SIAM Journal on Applied Mathematics 14, 641 - 664 (1966).

\bibitem{EY} L.C. Evans, Y. Yu, \emph{Various properties of solutions to the Infinity-Laplacian equation}, Communications in PDE 30:9, 1401 - 1428 (2005).
 
\bibitem{GNP} A. Garroni, V. Nesi, M. Ponsiglione, \emph{Dielectric breakdown: optimal bounds}, Proceedings of the Royal Society A  457, issue 2014 (2001).


\bibitem{J} R. Jensen, \emph{Uniqueness of Lipschitz extensions minimizing the sup-norm of the gradient}, Arch. Rational Mech. Analysis 123 (1993),
51-74.

\bibitem{K1} N. Katzourakis, \emph{$L^{\infty}$ Variational Problems for Maps and the Aronsson PDE System}, J. Differential Equations 253 (7), 2123-2139 (2012).

\bibitem{K1/2}  N. Katzourakis, \emph{Explicit $2D$ $\infty$-Harmonic Maps whose Interfaces have Junctions and Corners}, Comptes Rendus Acad. Sci. Paris, Ser.I, 351, 677 - 680 (2013).

\bibitem{K2} N. Katzourakis,  \emph{$\infty$-Minimal Submanifolds}, Proceedings of the AMS 142, 2797-2811 (2014).

\bibitem{K3}  N. Katzourakis, \emph{On the Structure of $\infty$-Harmonic Maps}, Communications in PDE 39:11, 2091 - 2124 (2014).

\bibitem{K4} N. Katzourakis, \emph{An Introduction to Viscosity Solutions for Fully Nonlinear PDE with Applications to Calculus of Variations in $L^\infty$}, Springer Briefs in Mathematics, 2015, DOI 10.1007/978-3-319-12829-0.

\bibitem{K5}  N. Katzourakis, \emph{Generalised solutions for fully nonlinear PDE systems and existence-uniqueness theorems}, Journal of Differential Equations 23, 641 - 686 (2017).

\bibitem{K6} N. Katzourakis,  \emph{Absolutely minimising generalised solutions to the equations of vectorial Calculus of Variations in $L^\infty$}, Calculus of Variations and PDE 56 (1), 1 - 25 (2017) (DOI: 10.1007/s00526-016-1099-z).

\bibitem{K7}  N. Katzourakis, \emph{A New Characterisation of $\infty$-Harmonic and $p$-Harmonic Mappings via Affine Variations in $L^\infty$}, Electronic Journal of Differential Equations 2017:29, 1 - 19 (2017).

\bibitem{KM} N. Katzourakis, R. Moser, \emph{Existence, Uniqueness and Structure of Second Order Absolute Minimisers}, Archives for Rational Mechanics and Analysis, published online 06/09/2018, DOI: 10.1007/s00205-018-1305-6.

\bibitem{KS}  N. Katzourakis, G. Shaw, \emph{Counterexamples in Calculus of Variations in $L^\infty$ through the vectorial Eikonal equation}, Comptes Rendus Mathematique Ser.\ I 356:5 (2018), 498-502, https://doi.org/10.1016/j.crma.2018.04.010.

\bibitem{KP} N. Katzourakis, T. Pryer, \emph{$2$nd order $L^\infty$ variational problems and the $\infty$-Polylaplacian}, Advances in Calculus of Variations, Published Online: 27-01-2018, DOI: https://doi.org/ 10.1515/acv-2016-0052.

\bibitem{MWZ} Q. Miao, C. Wang, Y. Zhou, \emph{Uniqueness of Absolute Minimizers for $L^\infty$-Functionals Involving Hamiltonians $H(x,p)$}, Archive for Rational Mechanics and Analysis 223 (1), 141-198 (2017).

\bibitem{PP} G. Papamikos, T. Pryer, \emph{A Lie symmetry analysis and explicit solutions of the two-dimensional $\infty$-Polylaplacian}, Studies in Applied Mathematics, Online 17 September 2018, https://doi.org/10.1111/sapm.12232.

\bibitem{P} F. Prinari, \emph{On the lower semicontinuity
and approximation of $L^\infty$-functionals}, NoDEA 22, 1591 - 1605 (2015).

\bibitem{RZ} A.N. Ribeiro, E. Zappale, \emph{Existence of minimisers for nonlevel convex functionals}, SIAM J. Control Opt., Vol. 52, No. 5,  (2014) 3341 - 3370.

\bibitem{SS} S. Sheffield, C.K. Smart, \emph{Vector Valued Optimal Lipschitz Extensions}, Comm. Pure Appl. Math. 65(1), 128 - 154 (2012).

\bibitem{Y} Y. Yu, \emph{Viscosity solutions of AronssonÕs equations}, Arch. Ration. Mech. Anal. 182, 153Ð180 (2006).

\end{thebibliography}

\end{document}